\documentclass[12pt]{article}%
\usepackage{amsfonts}
\usepackage{sw20bams}%
\usepackage{amsmath}%
\usepackage{amssymb}%
\usepackage{graphicx}
\providecommand{\U}[1]{\protect\rule{.1in}{.1in}}
\begin{document}

\title{Capturing, Ordering and Gaussianity in 2D}
\author{Steven Finch}
\date{January 19, 2016}
\maketitle

\begin{abstract}
We collect various facts related loosely to random Gaussian quadrilaterals in
the plane. \ For example, a side of a degenerate quadrilateral (one point
inside three others) has a density that is non-Rayleigh.

\end{abstract}

\footnotetext{Copyright \copyright \ 2016 by Steven R. Finch. All rights
reserved.}Let $A$, $B$, $C$ be independent random Gaussian points in
$\mathbb{R}^{2}$, all of which have mean vector zero and covariance matrix
identity. \ The triangle $ABC$ is said to \textbf{capture} the origin if
$(0,0)$ is contained in the convex hull of $A$, $B$, $C$. \ This event occurs
with probability $1/4$ and a highly attractive proof appears in
\cite{HS-CptrOrdG}. \ We offer an \textit{unattractive} proof of the same,
with the advantage that $(0,0)$ can be replaced by an arbitrary location
$(\xi,\eta)$, but the results are available only numerically. \ As far as is
known, a closed-form expression for probabilities does not exist.

We then review the variance of the median of three points in $\mathbb{R}^{1}$
and examine how this might be generalized to four points in $\mathbb{R}^{2}$.
\ The convex-hull peeling characterization for \textbf{order statistics}
(better: \textbf{data depth}) in the plane is usually said to be
distributionally intractible \cite{Sh-CptrOrdG, Ba-CptrOrdG, Al-CptrOrdG,
HR-CptrOrdG}. \ The reason will become rather clear here! \ Multivariate
integrals are written down, but no attempt is made to evaluate them.

Finally, we summarize what is known about random Gaussian quadrilaterals in
$\mathbb{R}^{2}$ and give some simulation-based outcomes. \ Our hope (as
always) is that someone else might be able to break the theoretical logjam and
provide a rigorous analysis supporting these.

\section{Location Capture}

Let $A=(a_{1},a_{2})$, $B=(b_{1},b_{2})$, $C=(c_{1},c_{2})$. \ The triangle
$ABC$ captures the location $(\xi,\eta)$ if and only if the three determinants
\cite{p1-CptrOrdG, p2-CptrOrdG, p3-CptrOrdG}%
\[%
\begin{array}
[c]{ccccc}%
\left\vert
\begin{array}
[c]{ccc}%
\xi & \eta & 1\\
b_{1} & b_{2} & 1\\
c_{1} & c_{2} & 1
\end{array}
\right\vert , &  & \left\vert
\begin{array}
[c]{ccc}%
a_{1} & a_{2} & 1\\
\xi & \eta & 1\\
c_{1} & c_{2} & 1
\end{array}
\right\vert , &  & \left\vert
\begin{array}
[c]{ccc}%
a_{1} & a_{2} & 1\\
b_{1} & b_{2} & 1\\
\xi & \eta & 1
\end{array}
\right\vert
\end{array}
\]
are all positive or all negative. \ By circular symmetry of the bivariate
normal distribution, we set $\eta=0$ without loss of generality (upon
rotation). \ Rather than work with $N(0,1)$ variables and $(\xi,0)$, we choose
$N(-\xi,1)$ variables and $(0,0)$ for simplicity's sake (upon translation).
\ Let us focus on the scenario%
\[%
\begin{array}
[c]{ccccc}%
b_{1}c_{2}-b_{2}c_{1}>0, &  & a_{1}b_{2}-a_{2}b_{1}>0, &  & a_{2}c_{1}%
-a_{1}c_{2}>0
\end{array}
\]
remembering to multiply by $2$ later. \ Within this, there are two cases:%
\[%
\begin{array}
[c]{ccccc}%
a_{1}>0, &  & b_{1}>0, &  & c_{1}<0;
\end{array}
\]%
\[%
\begin{array}
[c]{ccccc}%
a_{1}<0, &  & b_{1}<0, &  & c_{1}>0
\end{array}
\]
and we must remember to multiply by $3$ later. \ 

The first case gives rise to%
\[%
\begin{array}
[c]{ccccc}%
c_{2}>\dfrac{b_{2}c_{1}}{b_{1}}, &  & b_{2}>\dfrac{a_{2}b_{1}}{a_{1}}, &  &
c_{2}<\dfrac{a_{2}c_{1}}{a_{1}}%
\end{array}
\]
and inner integrals%
\begin{align*}
& \frac{1}{(2\pi)^{3/2}}%
{\displaystyle\int\limits_{-\infty}^{\infty}}
\;%
{\displaystyle\int\limits_{a_{2}b_{1}/a_{1}}^{\infty}}
\;%
{\displaystyle\int\limits_{b_{2}c_{1}/b_{1}}^{a_{2}c_{1}/a_{1}}}
\exp\left(  -\frac{a_{2}^{2}+b_{2}^{2}+c_{2}^{2}}{2}\right)  dc_{2}%
db_{2}da_{2}\\
& =\frac{1}{4\pi}%
{\displaystyle\int\limits_{-\infty}^{\infty}}
\;%
{\displaystyle\int\limits_{a_{2}b_{1}/a_{1}}^{\infty}}
\exp\left(  -\frac{a_{2}^{2}+b_{2}^{2}}{2}\right)  \left[  \operatorname{erf}%
\left(  \dfrac{a_{2}c_{1}}{\sqrt{2}a_{1}}\right)  -\operatorname{erf}\left(
\dfrac{b_{2}c_{1}}{\sqrt{2}b_{1}}\right)  \right]  db_{2}da_{2}\\
& =\frac{1}{4\sqrt{2\pi}}%
{\displaystyle\int\limits_{-\infty}^{\infty}}
\exp\left(  -\frac{a_{2}^{2}}{2}\right)  \left[  1+\operatorname{erf}\left(
\dfrac{a_{2}c_{1}}{\sqrt{2}a_{1}}\right)  +4T\left(  \dfrac{a_{2}c_{1}}{a_{1}%
},\frac{b_{1}}{c_{1}}\right)  \right]  da_{2}\\
& =\frac{1}{4}+\frac{1}{\sqrt{\pi}}%
{\displaystyle\int\limits_{-\infty}^{\infty}}
\exp\left(  -t^{2}\right)  T\left(  t\,z,w\right)  dt
\end{align*}
where $t=a_{2}/\sqrt{2}$, $z=\sqrt{2}c_{1}/a_{1}$, $w=b_{1}/c_{1}$ and
$T(h,k)$ is Owen's $T$ function \cite{Ow-CptrOrdG, FA-CptrOrdG}%
\[
T(h,k)=\frac{1}{2\pi}%
{\displaystyle\int\limits_{0}^{k}}
\frac{1}{1+s^{2}}\exp\left(  -\frac{h^{2}\left(  1+s^{2}\right)  }{2}\right)
ds.
\]
Integrating by parts:%
\[%
\begin{array}
[c]{lll}%
u=T\left(  t\,z,w\right)  , &  & dv=\dfrac{1}{\sqrt{\pi}}\exp\left(
-t^{2}\right)  dt,\\
du=-\dfrac{z}{2\sqrt{2\pi}}\exp\left(  -\dfrac{t^{2}z^{2}}{2}\right)
\operatorname{erf}\left(  \dfrac{t\,z\,w}{\sqrt{2}}\right)  dt, &  &
v=\dfrac{1}{2}\operatorname{erf}(t)
\end{array}
\]
we obtain%
\begin{align*}
& \frac{1}{4}+\dfrac{z}{4\sqrt{2\pi}}%
{\displaystyle\int\limits_{-\infty}^{\infty}}
\exp\left(  -\frac{t^{2}z^{2}}{2}\right)  \operatorname{erf}%
(t)\operatorname{erf}\left(  \dfrac{t\,z\,w}{\sqrt{2}}\right)  dt\\
& =\frac{1}{4}+\dfrac{p}{4\sqrt{\pi}}%
{\displaystyle\int\limits_{-\infty}^{\infty}}
\exp\left(  -p^{2}t^{2}\right)  \operatorname{erf}(t)\operatorname{erf}\left(
q\,t\right)  dt
\end{align*}
where $p=z/\sqrt{2}=c_{1}/a_{1}$ and $q=z\,w/\sqrt{2}=b_{1}/a_{1}$. No
antiderivative is possible, but the definite integral is known
\cite{Dk-CptrOrdG}, yielding%
\[
\frac{1}{4}+\dfrac{p}{4\sqrt{\pi}}\frac{2}{\sqrt{\pi}p}\arctan\left(  \frac
{q}{p\sqrt{1+p^{2}+q^{2}}}\right)  =\frac{1}{4}+\frac{1}{2\pi}\arctan\left(
\frac{a_{1}b_{1}}{c_{1}\sqrt{a_{1}^{2}+b_{1}^{2}+c_{1}^{2}}}\right)  .
\]

The second case gives rise to%
\[%
\begin{array}
[c]{ccccc}%
c_{2}<\dfrac{b_{2}c_{1}}{b_{1}}, &  & b_{2}<\dfrac{a_{2}b_{1}}{a_{1}}, &  &
c_{2}>\dfrac{a_{2}c_{1}}{a_{1}}%
\end{array}
\]
and inner integrals%
\begin{align*}
& \frac{1}{(2\pi)^{3/2}}%
{\displaystyle\int\limits_{-\infty}^{\infty}}
\;%
{\displaystyle\int\limits_{-\infty}^{a_{2}b_{1}/a_{1}}}
\;%
{\displaystyle\int\limits_{a_{2}c_{1}/a_{1}}^{b_{2}c_{1}/b_{1}}}
\exp\left(  -\frac{a_{2}^{2}+b_{2}^{2}+c_{2}^{2}}{2}\right)  dc_{2}%
db_{2}da_{2}\\
& =\frac{1}{4\pi}%
{\displaystyle\int\limits_{-\infty}^{\infty}}
\;%
{\displaystyle\int\limits_{-\infty}^{a_{2}b_{1}/a_{1}}}
\exp\left(  -\frac{a_{2}^{2}+b_{2}^{2}}{2}\right)  \left[  \operatorname{erf}%
\left(  \dfrac{b_{2}c_{1}}{\sqrt{2}b_{1}}\right)  -\operatorname{erf}\left(
\dfrac{a_{2}c_{1}}{\sqrt{2}a_{1}}\right)  \right]  db_{2}da_{2}\\
& =\frac{1}{4\sqrt{2\pi}}%
{\displaystyle\int\limits_{-\infty}^{\infty}}
\exp\left(  -\frac{a_{2}^{2}}{2}\right)  \left[  1-\operatorname{erf}\left(
\dfrac{a_{2}c_{1}}{\sqrt{2}a_{1}}\right)  +4T\left(  \dfrac{a_{2}c_{1}}{a_{1}%
},\frac{b_{1}}{c_{1}}\right)  \right]  da_{2}.
\end{align*}
The only difference with before is the sign in front of the error function;
the integral ultimately simplifies to%
\[
\frac{1}{4}-\frac{1}{2\pi}\arctan\left(  \frac{a_{1}b_{1}}{c_{1}\sqrt
{a_{1}^{2}+b_{1}^{2}+c_{1}^{2}}}\right)
\]
because $\left\vert a_{1}\right\vert =-a_{1}$ here.

Multiplying both expressions by $6$, the desired probability is%
\[
\frac{3}{(2\pi)^{5/2}}\left[  \varphi(\xi)+\psi(\xi)\right]  =\left\{
\begin{array}
[c]{lll}%
0.250000... &  & \text{if }\xi=0,\\
0.197171... &  & \text{if }\xi=1/2,\\
0.098289... &  & \text{if }\xi=1,\\
0.032455... &  & \text{if }\xi=3/2,\\
0.007626... &  & \text{if }\xi=2
\end{array}
\right.
\]
where%
\[
\varphi=%
{\displaystyle\int\limits_{0}^{\infty}}
{\displaystyle\int\limits_{0}^{\infty}}
{\displaystyle\int\limits_{-\infty}^{0}}
\exp\left(  -\tfrac{(a_{1}+\xi)^{2}+(b_{1}+\xi)^{2}+(c_{1}+\xi)^{2}}%
{2}\right)  \left[  \pi+2\arctan\left(  \tfrac{a_{1}b_{1}}{c_{1}\sqrt
{a_{1}^{2}+b_{1}^{2}+c_{1}^{2}}}\right)  \right]  dc_{1}db_{1}da_{1},
\]%
\[
\psi=%
{\displaystyle\int\limits_{-\infty}^{0}}
{\displaystyle\int\limits_{-\infty}^{0}}
{\displaystyle\int\limits_{0}^{\infty}}
\exp\left(  -\tfrac{(a_{1}+\xi)^{2}+(b_{1}+\xi)^{2}+(c_{1}+\xi)^{2}}%
{2}\right)  \left[  \pi-2\arctan\left(  \tfrac{a_{1}b_{1}}{c_{1}\sqrt
{a_{1}^{2}+b_{1}^{2}+c_{1}^{2}}}\right)  \right]  dc_{1}db_{1}da_{1}.
\]
Further symbolic integration does not appear to be possible.

In the preceding analysis, the location $(\xi,\eta)$ was assumed to be fixed.
\ Suppose instead that $X$ is a random Gaussian point in $\mathbb{R}^{2}$ --
independent of $A$, $B$, $C$ -- with mean vector zero and covariance identity.
\ The probability that $X$ falls in the convex hull of $A$, $B$, $C$ is
\cite{Ef-CptrOrdG, MC-CptrOrdG}%
\[
-\frac{1}{2}+\frac{3}{2}\frac{\operatorname{arcsec}(3)}{\pi}%
=0.0877398280459109052562833...=\frac{1-\theta}{4},
\]
the \textbf{expected probability content} of $ABC$. \ As the name suggests,
this is a mean (over all Gaussian triangles $ABC$ in the plane). \ The
corresponding variance is unknown, although recent progress has been made
\cite{Bu-CptrOrdG, BR-CptrOrdG}. \ We will see $\theta$ again later.

Similarly, a tetrahedron in $\mathbb{R}^{3}$ captures the origin iff four
$3\times3$ determinants are all positive or negative. \ One sample determinant
is%
\[
\left\vert
\begin{array}
[c]{cccc}%
a_{1} & a_{2} & a_{3} & 1\\
b_{1} & b_{2} & b_{3} & 1\\
c_{1} & c_{2} & c_{3} & 1\\
0 & 0 & 0 & 1
\end{array}
\right\vert =a_{1}b_{2}c_{3}+a_{2}b_{3}c_{1}+a_{3}b_{1}c_{2}-a_{1}b_{3}%
c_{2}-a_{2}b_{1}c_{3}-a_{3}b_{2}c_{1}
\]
and the algebraic challenge of working with four such inequalities is clear.
\ The expected probability content of a Gaussian tetrahedron in space is,
however, known to be \cite{Ef-CptrOrdG, MC-CptrOrdG}%
\[
-\frac{2}{5}+\frac{\operatorname{arcsec}(4)}{\pi}%
=0.0195693767448337562290498....
\]
This quantity is also called the \textbf{Gaussian volume}\ since it is a
distribution-dependent analog of ordinary Euclidean measure.

\section{Median Variance}

In one dimension, points $a$ and $b$ capture $x$ iff the two determinants%
\[%
\begin{array}
[c]{ccc}%
\left\vert
\begin{array}
[c]{cc}%
x & 1\\
b & 1
\end{array}
\right\vert , &  & \left\vert
\begin{array}
[c]{cc}%
a & 1\\
x & 1
\end{array}
\right\vert
\end{array}
\]
are both positive or both negative. \ Let us focus on the scenario $x-b<0$ and
$a-x<0$, that is, over the domain%
\[
\Omega_{a,b}=\left\{  x\in\mathbb{R}:a<x<b\right\}  .
\]
The desired variance is \cite{Jo-CptrOrdG}%
\[
\frac{6}{(2\pi)^{3/2}}%
{\displaystyle\int\limits_{\mathbb{R}}}
{\displaystyle\int\limits_{\mathbb{R}}}
{\displaystyle\int\limits_{\Omega_{a,b}}}
x^{2}\exp\left(  -\frac{a^{2}+b^{2}+x^{2}}{2}\right)  dx\,db\,da=1-\frac
{\sqrt{3}}{\pi}.
\]
A simulation suggests that the median is normally distributed; this is
\textit{false}, but the true density%
\[
\frac{6}{(2\pi)^{3/2}}%
{\displaystyle\int\limits_{-\infty}^{x}}
{\displaystyle\int\limits_{x}^{\infty}}
\exp\left(  -\frac{a^{2}+b^{2}+x^{2}}{2}\right)  db\,da=\frac{3}{2\sqrt{2\pi}%
}\exp\left(  -\frac{x^{2}}{2}\right)  \left[  1-\operatorname{erf}\left(
\frac{x}{\sqrt{2}}\right)  ^{2}\right]
\]
is nevertheless very close to that for $N(0,1-\sqrt{3}/\pi)$.

We have already stated requirements in two dimensions for points $A$, $B$, $C
$ to capture $X$. \ Define%
\[
\Omega_{A,B,C}=\left\{
\begin{array}
[c]{c}%
X\in\mathbb{R}^{2}:(b_{1}-a_{1})(x_{2}-a_{2})-(b_{2}-a_{2})(x_{1}-a_{1})>0,\\
\;\;\;\;\;\;\;\;\;\;\;\;\,(c_{1}-b_{1})(x_{2}-b_{2})-(c_{2}-b_{2})(x_{1}%
-b_{1})>0\\
\;\;\;\;\;\;\text{and}\;\,(a_{1}-c_{1})(x_{2}-c_{2})-(a_{2}-c_{2})(x_{1}%
-c_{1})>0
\end{array}
\right\}
\]
and let $\left\vert X\right\vert ^{2}=x_{1}^{2}+x_{2}^{2}$ for convenience.
\ The desired variance is%
\[
\sigma^{2}=\frac{8}{(2\pi)^{4}(1-\theta)}%
{\displaystyle\int\limits_{\mathbb{R}^{2}}}
{\displaystyle\int\limits_{\mathbb{R}^{2}}}
{\displaystyle\int\limits_{\mathbb{R}^{2}}}
{\displaystyle\int\limits_{\Omega_{A,B,C}}}
x_{1}^{2}\exp\left(  -\frac{\left\vert A\right\vert ^{2}+\left\vert
B\right\vert ^{2}+\left\vert C\right\vert ^{2}+\left\vert X\right\vert ^{2}%
}{2}\right)  dX\,dC\,dB\,dA.
\]
Simulation suggests normality again holds (approximately) and that $\sigma
^{2}\approx0.36$. \ Evaluating this eight-fold integral, or even removing just
$x_{2}$, does not seem feasible. \ The factor involving%
\[
1-\theta=-2+6\frac{\operatorname{arcsec}(3)}{\pi}%
=0.3509593121836436210251333...
\]
is necessary for reasons to be made clear shortly.

\section{Gaussian Quadrilaterals}

Let $A$, $B$, $C$, $D$ be independent random Gaussian points in $\mathbb{R}%
^{2}$, all of which have mean vector zero and covariance matrix identity.
\ The convex hull $ABCD$ of $A$, $B$, $C$, $D$ possesses either three or four
vertices; if the former is true, $ABCD$ is a \textbf{degenerate quadrilateral}%
; if the latter is true, $ABCD$ is a \textbf{quadrilateral}. \ (Non-degeneracy
rules out triangles.) \ The probability that the number of vertices is four is
\cite{Ef-CptrOrdG, MC-CptrOrdG}%
\[
\theta=3-6\frac{\operatorname{arcsec}(3)}{\pi}%
=0.6490406878163563789748666....
\]
Earlier occurrences of $\theta$ can now be explained. \ In Section 1,
$(1-\theta)/4$ is the probability that a specific point (say, $D$) is inside
the other three ($A$, $B$, $C$). \ By contrast, $1-\theta$ is the probability
that \textit{any} one of the four points is internal to the rest. In Section
2, given three points on a line, a median always exists (thus the denominator
is $1$). \ By contrast, given four points in the plane, we have implicitly
conditioned on degeneracy so that an inner point exists (thus the denominator
here is $1-\theta$).

With no constraints, the convex hull $ABCD$ has expected area $\sqrt{3}$ and
expected perimeter%
\[
(3+\theta)\sqrt{\pi}=6.4677562192310137839669010....
\]
It is interesting that $\sqrt{3}$ is twice the expected area of a planar
Gaussian triangle \cite{Fi-CptrOrdG, CG-CptrOrdG, Mi1-CptrOrdG, Mi2-CptrOrdG}%
,\ but the corresponding triangular perimeter $3\sqrt{\pi}$ is not so simply
related to $(3+\theta)\sqrt{\pi}$. \ No higher moments for either area of
$ABCD$ or perimeter of $ABCD$ are precisely known.

If the number of vertices of $ABCD$ is constrained to be three, then%
\[%
\begin{array}
[c]{ccc}%
\operatorname*{E}\left(  \text{side}\right)  \approx2.11, &  &
\operatorname*{E}\left(  \text{side}^{2}\right)  \approx5.32.
\end{array}
\]
(By \textquotedblleft side\textquotedblright\ is meant the distance between
two adjacent vertices on the exterior.) \ We deduce that the side density here
is non-Rayleigh as $\operatorname*{E}\left(  \text{side}\right)
^{2}/\operatorname*{E}\left(  \text{side}^{2}\right)  $ is not sufficiently
near to $\pi/4$. \ If the number of vertices of $ABCD$ is constrained to be
four, then%
\[%
\begin{array}
[c]{ccc}%
\operatorname*{E}\left(  \text{side}\right)  \approx1.45, &  &
\operatorname*{E}\left(  \text{side}^{2}\right)  \approx2.78.
\end{array}
\]
We are also interested in cross-correlations between sides. \ For the case of
quadrilaterals, how does the cross-correlation between sides sharing a vertex
differ from the cross-correlation between sides that are wholly disjoint?

Problems involving random quadrilaterals have occupied our attention for many
years \cite{Mi0-CptrOrdG} and remain essentially unsolved. \ Much more
relevant material can be found at \cite{Fn-CptrOrdG}, including experimental
computer runs that aided theoretical discussion here.


\begin{thebibliography}{99}                                                                                               %
\bibitem {HS-CptrOrdG}R. Howard and P. Sisson, Capturing the origin with
random points: generalizations of a Putnam problem, \textit{College Math. J.}
27 (1996) 186--192; MR1390366.

\bibitem {Sh-CptrOrdG}M. I. Shamos, Geometry and statistics: problems at the
interface, \textit{Algorithms and Complexity}, Proc. Carnegie-Mellon Univ.
sympos., ed. J. F. Traub, Academic Press, 1976, pp. 251--280; MR0431785 (55 \#4780).

\bibitem {Ba-CptrOrdG}V. Barnett, The ordering of multivariate data,
\textit{J. Royal Statist. Soc. Ser. A} 139 (1976) 318--355; MR0445726 (56 \#4060).

\bibitem {Al-CptrOrdG}G. Aloupis, Geometric measures of data depth,
\textit{Data Depth: Robust Multivariate Analysis, Computational Geometry and
Applications}, Proc. 2003 Rutgers Univ. workshop, ed. R. Y. Liu, R. Serfling
and D. L. Souvaine, Amer. Math. Soc., 2006, pp. 147--158; MR2343118.

\bibitem {HR-CptrOrdG}J. Hugg, E. Rafalin, K. Seyboth and D. Souvaine, An
experimental study of old and new depth measures, \textit{Proc. 2006 Workshop
on Algorithm Engineering and Experiments (ALENEX)}, Miami, ed. R. Raman and M.
F. Stallmann, SIAM, pp. 51--64.

\bibitem {p1-CptrOrdG}G. Herron, Point within a tetrahedron,
USENET\ comp.graphics.algorithms posting, 1994.

\bibitem {p2-CptrOrdG}Anonymous contributor, Determining if an arbitrary point
lies inside a triangle defined by three points, Mathematics Stack Exchange
posting, 2011.

\bibitem {p3-CptrOrdG}Anonymous contributor, How to check whether the point is
in the tetrahedron or not, Stack Overflow posting, 2014.

\bibitem {Ow-CptrOrdG}D. B. Owen, Tables for computing bivariate normal
probabilities, \textit{Annals Math. Statist.} 27 (1956) 1075--1090; MR0127562
(23 \#B607).

\bibitem {FA-CptrOrdG}H. A. Fayed and A. F. Atiya, An evaluation of the
integral of the product of the error function and the normal probability
density with application to the bivariate normal integral, \textit{Math.
Comp}. 83 (2014) 235--250; MR3120588.

\bibitem {Dk-CptrOrdG}A. Dieckmann, Table of Integrals,
http://www-elsa.physik.uni-bonn.de/\symbol{126}dieckman/IntegralsDefinite/DefInt.html.

\bibitem {Ef-CptrOrdG}B. Efron, The convex hull of a random set of points,
\textit{Biometrika} 52 (1965) 331--343; MR0207004 (34 \#6820).

\bibitem {MC-CptrOrdG}S. N. Majumdar, A. Comtet and J. Randon-Furling, Random
convex hulls and extreme value statistics, \textit{J. Stat. Phys.} 138 (2010)
955--1009; MR2601420 (2011c:62166).

\bibitem {Bu-CptrOrdG}C. Buchta, An identity relating moments of functionals
of convex hulls, \textit{Discrete Comput. Geom.} 33 (2005) 125--142; MR2105754 (2005j:60022).

\bibitem {BR-CptrOrdG}M. Beermann and M. Reitzner, Beyond the Efron-Buchta
identities: distributional results for Poisson polytopes, \textit{Discrete
Comput. Geom.} 53 (2015) 226--244; MR3293497.

\bibitem {Jo-CptrOrdG}H. L. Jones, Exact lower moments of order statistics in
small samples from a normal distribution, \textit{Annals Math. Statist}. 19
(1948) 270--273; MR0025121 (9,601d).

\bibitem {Fi-CptrOrdG}S. R. Finch, Random triangles, unpublished note (2010),
http://www.people.fas.harvard.edu/\symbol{126}sfinch/.

\bibitem {CG-CptrOrdG}P. Clifford and N. J. B. Green, Distances in Gaussian
point sets, \textit{Math. Proc. Cambridge Philos. Soc.} 97 (1985) 515--524;
MR0778687 (86i:62091).

\bibitem {Mi1-CptrOrdG}K. S. Miller, Complex Gaussian processes, \textit{SIAM
Rev.} 11 (1969) 544--567; MR0258109 (41 \#2756).

\bibitem {Mi2-CptrOrdG}K. S. Miller, \textit{Complex Stochastic Processes. An
Introduction to Theory and Application}, Addison-Wesley, 1974, pp. 86--100;
MR0368118 (51 \#4360).

\bibitem {Mi0-CptrOrdG}W. J. C. Miller, Problem 6543, \textit{The Educational
Times} 33 (1880) 311; 53 (1900) 225.

\bibitem {Fn-CptrOrdG}S. R. Finch, Simulations in R\ involving quadrilaterals,
http://www.people.fas.harvard.edu/\symbol{126}sfinch/csolve/rsimul.html.%

\begin{tabular}
[c]{lll}
& Steven Finch & \\
& Dept. of Statistics & \\
& Harvard University & \\
& Cambridge, MA, USA & \\
& \textit{steven\_finch@harvard.edu} &
\end{tabular}

\end{thebibliography}
\end{document}